\documentclass[12pt]{article}

\usepackage{amsmath, epsfig, cite}
\usepackage{amssymb}
\usepackage{amsfonts}
\usepackage{latexsym}
\usepackage{graphicx}

\newtheorem{thm}{Theorem}[section]
\newtheorem{prop}[thm]{Proposition}

\newtheorem{cor}[thm]{Corollary}

\newtheorem{lem}[thm]{Lemma}

\newtheorem{rem}[thm]{Remark}

\newcommand{\qed}{{\hfill\rule{4pt}{7pt}}}

\numberwithin{equation}{section}

\makeatletter \@addtoreset{equation}{section} \makeatother

\begin{document}

\rule{0cm}{1cm}

\begin{center}
{\Large\bf The Cauchy Operator  for

 Basic Hypergeometric Series}
\end{center}
\vskip 2mm \centerline{Vincent Y. B. Chen$^1$ and Nancy S. S. Gu$^2$
}

\begin{center}
Center for Combinatorics, LPMC\\
Nankai University, Tianjin 300071\\
People's Republic of China\\

\vskip 2mm
 Email:  $^1$ybchen@mail.nankai.edu.cn,
 $^2$gu@nankai.edu.cn
\end{center}

\begin{center}
{\bf Abstract}
\end{center}

{\small We introduce the Cauchy augmentation operator for basic hypergeometric
series. Heine's ${}_2\phi_1$ transformation formula and Sears'
${}_3\phi_2$ transformation formula can be easily obtained by the
symmetric property of some parameters in operator identities. The
Cauchy  operator involves two parameters, and it
 can be considered as a generalization of the operator $T(bD_q)$.
Using this operator, we obtain extensions of the Askey-Wilson
integral, the Askey-Roy integral, Sears' two-term summation formula,
as well as the $q$-analogues of Barnes' lemmas. Finally, we find
that the Cauchy operator is also suitable for the study of the
bivariate Rogers-Szeg\"o polynomials, or the continuous big
$q$-Hermite polynomials.}

\vskip 5mm



\section{Introduction}

In an attempt to find efficient $q$-shift operators to deal with
basic hypergeometric series identities in the framework of the
$q$-umbral calculus \cite{Andrews, Goldman-Rota},  Chen and Liu
\cite{Chen1,Chen2} introduced two $q$-exponential operators  for
deriving identities from their special cases. This method is called
parameter augmentation.  In this paper, we continue the study of
parameter augmentation by defining a new operator called the Cauchy
augmentation operator which is suitable for certain transformation
and integral formulas.

Recall that Chen and Liu \cite{Chen1} introduced the augmentation
operator
\begin{equation}\label{TE}
T(bD_q)=\sum_{n=0}^{\infty}\frac{(bD_q)^n}{(q;q)_n}
\end{equation}
as the basis of parameter augmentation which serves as a method for
proving $q$-summation and integral formulas from special cases for
which  some parameters are set to zero.

The main idea of this paper is to introduce the Cauchy augmentation
operator, or simply the Cauchy operator,
\begin{equation}\label{gTE}
T(a,b;D_q)=\sum_{n=0}^{\infty}\frac{(a;q)_n}{(q;q)_n}(bD_q)^n,
\end{equation}
which is reminiscent of the Cauchy
 $q$-binomial
theorem \cite[Appendix II.3]{Gasper-Rahman}
\begin{equation}\label{Cauchy}
\sum_{n=0}^{\infty}\frac{(a;q)_n}{(q;q)_n}z^n
=\frac{(az;q)_\infty}{(z;q)_\infty},\ \ |z|<1.
\end{equation}
For the same reason, the operator $T(aD_q)$ should be named the
Euler operator in view of Euler's identity\cite[Appendix
II.1]{Gasper-Rahman}
\begin{equation}\label{Euler}
\sum_{n=0}^{\infty}\frac{z^n}{(q;q)_n} =\frac{1}{(z;q)_\infty},\ \ \
\ \quad |z|<1.
\end{equation}

 Compared with $T(bD_q)$, the Cauchy
 operator \eqref{gTE} involves two parameters. Clearly,
the operator $T(bD_q)$ can be considered as a special case of the
Cauchy operator \eqref{gTE} for $a=0$. In order to utilize the
Cauchy  operator to basic hypergeometric series, several operator
identities are deduced in Section \ref{se Basic Properties}. As to
the applications of the Cauchy operator, we show that many classical
results on basic hypergeometric series easily fall into this
framework. Heine's ${}_2\phi_1$ transformation formula
\cite[Appendix III.2]{Gasper-Rahman} and Sears' ${}_3\phi_2$
transformation formula \cite[Appendix III.9]{Gasper-Rahman} can be
easily obtained by the symmetric property of some parameters in two
operator identities for the Cauchy operator.

In Section \ref{se Askey-Wilson Integral} and Section \ref{se
Askey-Roy Integral}, we use the Cauchy  operator to generalize the
Askey-Wilson integral and the Askey-Roy integral. In
\cite{IsmailStantonViennot}, Ismail, Stanton, and Viennot derived an
integral named the Ismail-Stanton-Viennot integral which took the
Askey-Wilson integral as a special case. It is easy to see that our
extension of the Askey-Wilson integral is also an extension of the
Ismail-Stanton-Viennot integral. In \cite{Gasper}, Gasper discovered
an integral which was a generalization of the Askey-Roy integral. We
observe that Gasper's formula is a special case of the formula
obtained by applying the Cauchy  operator directly to the Askey-Roy
integral. Furthermore, we find that the Cauchy  operator can be
applied to Gasper's formula to derive a further extension of the
Askey-Roy integral.

In Section \ref{se Bivariate Rogers-Szego Polynomials}, we present
that the Cauchy  operator is suitable for the study of bivariate
Rogers-Szeg\"o polynomials. It can be used to derive the
corresponding Mehler's and the Rogers formulas for the bivariate
Rogers-Szeg\"o polynomials, which can be stated in the equivalent
forms in terms of the continuous big $q$-Hermite polynomials.
Mehler's formula in this case turns out to be a special case of the
nonsymmetric Poisson kernel formula for the continuous big
$q$-Hermite polynomials due to Askey, Rahman, and Suslov
\cite{Askey-Rahman-Suslov}. Finally, in Section \ref{se Sears'
Formula} and Section \ref{se q-Barnes' Lemmas}, we employ the Cauchy
operator to deduce extensions of Sears' two-term summation formula
\cite[Eq. (2.10.18)]{Gasper-Rahman} and the $q$-analogues of Barnes'
lemmas \cite[Eqs. (4.4.3), (4.4.6)]{Gasper-Rahman}.

As usual, we follow the notation and terminology in
\cite{Gasper-Rahman}. For $|q|<1$, the $q$-shifted factorial is
defined by
$$(a;q)_\infty=
\prod_{k=0}^{\infty}(1-aq^k) \text{\ \  and \ \  }(a;q)_n
=\frac{(a;q)_\infty}{(aq^n;q)_\infty}, \text{ for } n\in
\mathbb{Z}.$$

For convenience, we shall adopt the following notation for multiple
$q$-shifted factorials:
$$(a_1,a_2,\ldots,a_m;q)_n=(a_1;q)_n(a_2;q)_n\cdots(a_m;q)_n,$$
where $n$ is an integer or infinity.

The $q$-binomial coefficients, or the Gauss coefficients, are given
by
\begin{equation}
{n \brack k}=\frac{(q;q)_n}{(q;q)_k(q;q)_{n-k}}.
\end{equation}

The (unilateral) basic hypergeometric series $_{r}\phi_s$ is defined
by
\begin{equation}
_{r}\phi_s\left[\begin{array}{cccccc} a_1,&a_2,&\ldots,&a_r\\
b_1,&b_2,&\ldots,&b_s
\end{array};q,z\right]=\sum_{k=0}^{\infty}
\frac{(a_1,a_2,\ldots,a_r;q)_k}{(q,b_1,b_2,\ldots,b_s;q)_k}
\left[(-1)^kq^{k\choose2}\right]^{1+s-r}z^k.
\end{equation}

\section{Basic Properties}\label{se Basic Properties}

In this section, we give some basic identities involving the Cauchy
 operator $T(a,b; D_q)$ and demonstrate that Heine's ${}_2\phi_1$
 transformation formula and Sears' ${}_3\phi_2$
transformation formula are implied in the symmetric property of some
parameters in two operator identities.

We recall that the $q$-difference operator, or Euler derivative, is
defined by
\begin{equation}\label{dq}
D_q\{f(a)\}=\frac{f(a)-f(aq)}{a},
\end{equation}
and the Leibniz rule for $D_q$ is referred to
 the following
identity
\begin{equation}\label{Leibniz}
D_q^n\{f(a)g(a)\}=\sum_{k=0}^n q^{k(k-n)}{n \brack
k}D_q^k\{f(a)\}D_q^{n-k}\{g(aq^k)\}.
\end{equation}
The following relations are easily verified.

\begin{prop}\label{DK1}Let $k$ be a nonnegative integer. Then we
have
\begin{eqnarray*}
D_q^k\left\{\frac{1}{(at;q)_\infty}\right\}
&=&\frac{t^k}{(at;q)_\infty},\\[5pt]
D_q^k\left\{(at;q)_\infty\right\}
& = & (-t)^kq^{k\choose2}(atq^k;q)_\infty,\\[5pt]
D_q^k\left\{\frac{(av;q)_\infty}{(at;q)_\infty}\right\}
& = & t^k(v/t;q)_k\frac{(avq^k;q)_\infty}{(at;q)_\infty}.\\
\end{eqnarray*}
\end{prop}

 Now, we are
ready to give some  basic identities for the Cauchy operator
$T(a,b;D_q)$. We assume that $T(a,b;D_q)$ acts on the parameter $c$.
The following   identity  is an easy consequence of the Cauchy
$q$-binomial theorem \eqref{Cauchy}.

\begin{thm}\label{gtef1}  We have
\begin{equation}\label{tabc}
T(a,b;D_q)\left\{\frac{1}{(ct;q)_\infty}\right\}=\frac{(ab\,t;q)_\infty}{(b\,t,ct;q)_\infty},
\end{equation}
provided $|b\,t|<1$.
\end{thm}

\begin{pf} By Proposition \ref{DK1}, the left hand side of \eqref{tabc} equals
\[ \sum_{n=0}^{\infty}\frac{(a;q)_nb^n}{(q;q)_n}
D_q^n\left\{\frac{1}{(ct;q)_\infty}\right\} =
 {1\over (ct;q)_\infty} \sum_{n=0}^{\infty}
 \frac{(a;q)_n(b\,t)^n}{(q;q)_n},\]
 which simplifies to the right hand side of \eqref{tabc}
 by the Cauchy $q$-binomial theorem \eqref{Cauchy}.\qed
\end{pf}

\begin{thm}\label{gtef2} We have
\begin{equation}\label{22}
T(a,b;D_q)\left\{\frac{1}{(cs,ct;q)_\infty}\right\}
=\frac{(ab\,t;q)_\infty}{(b\,t,cs,ct;q)_\infty}
\,{}_{2}\phi_1\left[\begin{array}{cc}
a,&ct\\
&ab\,t
\end{array};q,bs\right],
\end{equation}
provided $\max\{|bs|,|b\,t|\}<1$.
\end{thm}

\begin{pf} In view of the Leibniz formula for $D_q^n$, the left
hand side of \eqref{22} can be expanded as follows
\begin{eqnarray*}
\lefteqn{\sum_{n=0}^{\infty}\frac{(a;q)_nb^n}{(q;q)_n}\sum_{k=0}^{n}q^{k(k-n)}
{n\brack k}D_q^k\left\{\frac{1}{(cs;q)_\infty}\right\}D_q^{n-k}
\left\{\frac{1}{(ctq^k;q)_\infty}\right\}\nonumber}\\[6pt]
&=&\sum_{n=0}^{\infty}\frac{(a;q)_nb^n}{(q;q)_n}\sum_{k=0}^{n}q^{k(k-n)}{n\brack
k}
\frac{s^k}{(cs;q)_\infty}\frac{(tq^k)^{n-k}}{(ctq^k;q)_\infty}
\nonumber\\[6pt]
&=&\frac{1}{(cs,ct;q)_\infty}\sum_{k=0}^{\infty}\frac{(ct;q)_k(bs)^k}{(q;q)_k}
\sum_{n=k}^{\infty}\frac{(a;q)_n(b\,t)^{n-k}}{(q;q)_{n-k}}
\nonumber\\[6pt]
&=&\frac{1}{(cs,ct;q)_\infty}\sum_{k=0}^{\infty}\frac{(a,ct;q)_k(bs)^k}{(q;q)_k}
\sum_{n=0}^{\infty}\frac{(aq^k;q)_n(b\,t)^{n}}{(q;q)_{n}}
\quad\quad\quad\quad\quad\quad\nonumber\\[6pt]
&=&\frac{(ab\,t;q)_\infty}{(b\,t,cs,ct;q)_\infty}
{}_{2}\phi_1\left[\begin{array}{cc}
a,&ct\\
&ab\,t
\end{array};q,bs\right],
\end{eqnarray*}
as desired. \qed
\end{pf}

Notice that when $a=0$, the ${}_2\phi_1$ series
 on the right hand
side of \eqref{22} can be summed by employing the Cauchy
$q$-binomial theorem \eqref{Cauchy}. In this case \eqref{22} reduces
to
\begin{equation}
T(bD_q)\left\{\frac{1}{(cs,ct;q)_\infty}\right\}=
\frac{(bcst;q)_\infty}{(bs,b\,t,cs,ct;q)_\infty},\quad
\quad|bs|,|b\,t|<1,
\end{equation}
which was derived by Chen and Liu in \cite{Chen1}.

As an immediate consequence of the above theorem, we see that
Heine's ${}_2\phi_1$ transformation formula \cite[Appendix
III.2]{Gasper-Rahman} is really about the symmetry in $s$ and $t$
while applying the operator $T(a,b;q)$.

\begin{cor}[Heine's transformation]\label{Heine}We have
\begin{equation}\label{sHeine}
{}_{2}\phi_1\left[\begin{array}{cc}
a,&b\\
&c
\end{array};q,z\right]=\frac{(c/b,bz;q)_\infty}
{(c,z;q)_\infty}{}_{2}\phi_1\left[\begin{array}{cc}
abz/c,&b\\
&bz
\end{array};q,\frac{c}{b}\right],
\end{equation}
where $\max\{|z|, |c/b|\}<1$.
\end{cor}

\begin{pf} The symmetry in $s$ and $t$ on the
left hand side of \eqref{22} implies that
\begin{equation}\label{21}
\frac{(ab\,t;q)_\infty}{(b\,t,cs,ct;q)_\infty}
{}_{2}\phi_1\left[\begin{array}{cc}
a,&ct\\
&ab\,t
\end{array};q,bs\right]=\frac{(abs;q)_\infty}{(bs,ct,cs;q)_\infty}
{}_{2}\phi_1\left[\begin{array}{cc}
a,&cs\\
&abs
\end{array};q,b\,t\right],
\end{equation}
where $\max\{|bs|,|b\,t|\}<1$.

Replacing $a, b, c, s, t$ by $b,a,a^2b/c, z/a, c/ab$ in \eqref{21},
respectively, we may easily express the above identity in the form
of \eqref{sHeine}.\qed
\end{pf}

\begin{rem}
 A closer look at the proof of Theorem \ref{gtef2} reveals that
the essence of Heine's transformation lies in the symmetry of $f$
and $g$ in Leibniz's formula \eqref{Leibniz}.
\end{rem}

 We should note that we must be cautious about the convergence
conditions while utilizing the Cauchy  operator. In general, it
would be safe to apply the Cauchy  operator if the resulting series
is convergent. However, it is possible that from a convergent series
one may obtain a divergent series after employing the Cauchy
 operator. For example, let us consider Corollary
\ref{Heine}. The resulting series \eqref{21} can be obtained  by
applying the Cauchy  operator $T(a,b;D_q)$ to $1/(cs,ct;q)_\infty$
which is convergent for all $t$.  However, the resulting series on
the left hand side of \eqref{21} is not convergent for $|t|>1/|b|$.

 Combining
Theorem \ref{gtef1} and the Leibniz rule  \eqref{Leibniz}, we obtain
the following identity which implies Theorem \ref{gtef2} by setting
$v=0$.  Sears' ${}_3\phi_2$ transformation formula \cite[Appendix
III.9]{Gasper-Rahman} is also a consequence of Theorem \ref{gtef3}.

\begin{thm}We have\label{gtef3}
\begin{eqnarray}\label{tab}
T(a,b;D_q)\left\{\frac{(cv;q)_\infty}{(cs,ct;q)_\infty}\right\}
=\frac{(abs,cv;q)_\infty}{(bs,cs,ct;q)_\infty}{}_{3}\phi_2\left[\begin{array}{ccc}
a,&cs,&v/t\\
&abs,&cv
\end{array};q,b\,t\right],
\end{eqnarray}
provided $\max\{|bs|,|b\,t|\}<1$.
\end{thm}

\begin{pf} In light of Leibniz's formula, the left hand side of
\eqref{tab} equals
\begin{eqnarray}
\lefteqn{\sum_{n=0}^{\infty}\frac{(a;q)_nb^n}{(q;q)_n}
D_q^n\left\{\frac{(cv;q)_\infty}{(cs,ct;q)_\infty}\right\}\nonumber}
\\[6pt]
&=&\sum_{n=0}^{\infty}\frac{(a;q)_nb^n}{(q;q)_n}
\sum_{k=0}^{n}q^{k(k-n)}
{n\brack k}D_q^k\left\{\frac{(cv;q)_\infty}{(ct;q)_\infty}\right\}
D_q^{n-k}\left\{\frac{1}{(csq^k;q)_\infty}\right\}\nonumber\\[6pt]
&=&\sum_{n=0}^{\infty}\frac{(a;q)_nb^n}{(q;q)_n}
\sum_{k=0}^{n}q^{k(k-n)}{n\brack
k} \frac{t^k(v/t;q)_k(cvq^k;q)_\infty}{(ct;q)_\infty}
D_q^{n-k}\left\{\frac{1}{(csq^k;q)_\infty}\right\}\nonumber\\[6pt]
&=&\sum_{k=0}^{\infty}\frac{(v/t;q)_k(cvq^k;q)_\infty{t}^k}
{(q;q)_k(ct;q)_\infty}
\sum_{n=k}^{\infty}\frac{b^nq^{k(k-n)}(a;q)_n}{(q;q)_{n-k}}
D_q^{n-k}\left\{\frac{1}{(csq^k;q)_\infty}\right\}\nonumber\\[6pt]
&=&\sum_{k=0}^{\infty}\frac{(a,v/t;q)_k(cvq^k;q)_\infty(b\,t)^k}
{(q;q)_k(ct;q)_\infty}
\sum_{n=0}^{\infty}\frac{(bq^{-k})^n(aq^k;q)_n}{(q;q)_n}
D_q^{n}\left\{\frac{1}{(csq^k;q)_\infty}\right\}\nonumber\\[6pt]
&=&\sum_{k=0}^{\infty}\frac{(a,v/t;q)_k(cvq^k;q)_\infty(b\,t)^k}{(q;q)_k(ct;q)_\infty}
T(aq^k,bq^{-k};D_q)\left\{\frac{1}{(csq^k;q)_\infty}\right\}
\nonumber.
\end{eqnarray}
By Theorem \ref{gtef1}, the above sum equals
\begin{eqnarray}
\lefteqn{\sum_{k=0}^{\infty}\frac{(a,v/t;q)_k(cvq^k;q)_\infty(b\,t)^k}{(q;q)_k(ct;q)_\infty}
\frac{(absq^k;q)_\infty}{(bs,csq^k;q)_\infty}\nonumber}\\[6pt]
&=&\frac{(cv;q)_\infty}{(cs,ct;q)_\infty}\sum_{k=0}^{\infty}
\frac{(a,cs,v/t;q)_k(b\,t)^k}{(q,cv;q)_k}
\frac{(absq^k;q)_\infty}{(bs;q)_\infty}\nonumber\\[6pt]
&=&\frac{(abs,cv;q)_\infty}
{(bs,cs,ct;q)_\infty}{}_{3}\phi_2\left[\begin{array}{ccc}
a,&cs,&v/t\\
&abs,&cv \label{33}
\end{array};q,b\,t\right],
\end{eqnarray}
as desired. \qed
\end{pf}

\begin{cor}[Sears' transformation]We have
\begin{equation}
{}_{3}\phi_2\left[\begin{array}{ccc}
a,&b,&c\\
&d,&e
\end{array};q,\frac{de}{abc}\right]=
\frac{(e/a,de/bc;q)_\infty}{(e,de/abc;q)_\infty}
{}_{3}\phi_2\left[\begin{array}{ccc}
a,&d/b,&d/c\\
&d,&de/bc
\end{array};q,\frac{e}{a}\right],\label{hall}
\end{equation}
where $\max\{|de/abc|,|e/a|\}<1$.
\end{cor}

\begin{pf}
Based on the symmetric property of the parameters $s$ and $t$ on the
left hand side of \eqref{tab}, we find that
\begin{equation*}
\frac{(abs,cv;q)_\infty}{(bs,cs,ct;q)_\infty}{}_{3}\phi_2\left[\begin{array}{ccc}
a,&cs,&v/t\\
&abs,&cv
\end{array};q,b\,t\right]=\frac{(ab\,t,cv;q)_\infty}{(b\,t,ct,cs;q)_\infty}
{}_{3}\phi_2\left[\begin{array}{ccc}
a,&ct,&v/s\\
&ab\,t,&cv
\end{array};q,bs\right],
\end{equation*}
where $\max\{|bs|,|b\,t|\}<1$.

Making the substitutions $c \rightarrow ab^2/e$, $v \rightarrow
de/ab^2$, $s \rightarrow e/ab$, and $t\rightarrow de/ab^2c$, we get
the desired formula.\qed
\end{pf}

We see that the essence of Sears' transformation also lies
in the symmetry of $s$ and $t$ in the application of Leibniz
rule.

\section{An Extension of the Askey-Wilson Integral}
\label{se Askey-Wilson Integral}

The Askey-Wilson integral \cite{AskeyWilson} is a significant
extension of the beta integral. Chen and Liu \cite{Chen1} presented
a treatment of the Askey-Wilson integral via parameter augmentation.
They first got the usual Askey-Wilson integral with one parameter by
the orthogonality relation obtained from the Cauchy $q$-binomial
theorem \eqref{Cauchy} and the Jacobi triple product identity
\cite[Appendix II.28]{Gasper-Rahman}, and then they applied the
operator $T(bD_q)$ three times to deduce the Askey-Wilson integral
involving four parameters \cite{Askey,IsmailStantonViennot,
IsmailStanton,Kalnins,Rahman,WilfZeilberger}
\begin{eqnarray}
&&\int_{0}^{\pi}\frac{(e^{2i\theta},e^{-2i\theta};q)_\infty d\theta}
{(ae^{i\theta},ae^{-i\theta},be^{i\theta},be^{-i\theta},
ce^{i\theta},ce^{-i\theta},de^{i\theta},de^{-i\theta};q)_\infty}\nonumber\\[6pt]
&&=\frac{2\pi(abcd;q)_\infty}{(q,ab,ac,ad,bc,bd,cd;q)_\infty},\label{Askeywilson}
\end{eqnarray}
where $\max\{|a|,|b|,|c|,|d|\}<1$.

In this section, we derive an extension of the Askey-Wilson integral
\eqref{Askeywilson}   which contains the following
Ismail-Stanton-Viennot's integral \cite{IsmailStantonViennot} as a
special case:
\begin{eqnarray}\label{ISV}
\lefteqn{\int_{0}^{\pi}\frac{(e^{2i\theta},e^{-2i\theta};q)_\infty
d\theta}
{(ae^{i\theta},ae^{-i\theta},be^{i\theta},be^{-i\theta},ce^{i\theta},ce^{-i\theta},
de^{i\theta},de^{-i\theta},ge^{i\theta},ge^{-i\theta};q)_\infty}}\nonumber\\[6pt]
&=&\frac{2\pi(abcg,abcd;q)_\infty}{(q,ab,ac,ad,ag,bc,bd,bg,cd,cg;q)_\infty}
{}_{3}\phi_2\left[\begin{array}{ccc}
ab,&ac,&bc\\
&abcg,&abcd
\end{array};q,dg\right],
\end{eqnarray}
where $\max\{|a|,|b|,|c|,|d|,|g|\}<1$.

\begin{thm}[Extension of the Askey-Wilson integral]We have
\label{ThmAskeyWilson}
\begin{eqnarray}\label{exAskeywilson}
\lefteqn{\int_{0}^{\pi}\frac{(e^{2i\theta},e^{-2i\theta},fge^{i\theta};q)_\infty}
{(ae^{i\theta},ae^{-i\theta},be^{i\theta},be^{-i\theta},ce^{i\theta},ce^{-i\theta},
de^{i\theta},de^{-i\theta},ge^{i\theta};q)_\infty}}\nonumber\\[6pt]
&&\quad\quad\quad\times{}_{3}\phi_2\left[\begin{array}{ccc}
f,&ae^{i\theta},&be^{i\theta}\\
&fge^{i\theta},&ab
\end{array};q,ge^{-i\theta}\right]d\theta\nonumber\\[6pt]
&=&\frac{2\pi(cfg,abcd;q)_\infty}{(q,ab,ac,ad,bc,bd,cd,cg;q)_\infty}
{}_{3}\phi_2\left[\begin{array}{ccc}
f,&ac,&bc\\
&cfg,&abcd
\end{array};q,dg\right],
\end{eqnarray}
where $\max\{|a|,|b|,|c|,|d|,|g|\}<1$.
\end{thm}

\begin{pf}
The Askey-Wilson integral \eqref{Askeywilson} can be written as
\begin{eqnarray}\label{integral}
&&\int_{0}^{\pi}\frac{(e^{2i\theta},e^{-2i\theta};q)_\infty}
{(be^{i\theta},be^{-i\theta},ce^{i\theta},ce^{-i\theta},de^{i\theta},de^{-i\theta};q)_\infty}
\frac{(ab;q)_\infty}{(ae^{i\theta},ae^{-i\theta};q)_\infty}d\theta\nonumber\\[6pt]
&&\qquad
=\frac{2\pi}{(q,bc,bd,cd;q)_\infty}\frac{(abcd;q)_\infty}{(ac,ad;q)_\infty}.
\end{eqnarray}

Before applying the Cauchy operator to an integral, it is necessary
to show that the Cauchy operator commutes with the integral. This
fact is implicit in the literature. Since this commutation relation
depends on some technical conditions in connection with the
integrands, here we present a complete proof.

First, it can be easily verified that the $q$-difference operator
$D_q$ commutes with the integral. By the definition of $D_q$
\eqref{dq}, it is clear  that
\begin{equation}
D_q\left\{\int_C f(\theta,a){d}\theta\right\} = \int_C D_q
\left\{f(\theta,a)\right\}{d}\theta.
\end{equation}
Consequently, the operator $D_q^n$ commutes with the integral.
Given a Cauchy operator
$T(f,g;D_q)$, we proceed to prove that it commutes with
the integral. From the well-known fact
that, for a sequence of continuous functions $u_n(\theta)$ on a
curve $C$, the sum commutes with the integral in
\[ \sum_{n=0}^{\infty}\int_{C}u_n(\theta){d}\theta\]
 provided that
$\sum_{n=0}^{\infty}u_n(\theta)$ is uniformly convergent.
It is
sufficient to check the convergence condition for the continuity is
obvious. This can be done with the aid of the Weierstrass M-Test
\cite{Arfken}. Using the Cauchy operator
$T(f,g;D_q)$ to the left hand side of \eqref{integral}, we find that
\begin{eqnarray}
\lefteqn{T(f,g;D_q)\left\{\int_{0}^{\pi}\frac{(e^{2i\theta},e^{-2i\theta};q)_\infty}
{(be^{i\theta},be^{-i\theta},ce^{i\theta},ce^{-i\theta},de^{i\theta},de^{-i\theta};q)_\infty}
\frac{(ab;q)_\infty}{(ae^{i\theta},ae^{-i\theta};q)_\infty}{d}\theta\right\}
}\label{44}\\[6pt]
&=& \sum_{n=0}^{\infty}\frac{(f;q)_{n}}{(q;q)_{n}}(gD_q)^n
\int_{0}^{\pi}\frac{(e^{2i\theta},e^{-2i\theta};q)_\infty}
{(be^{i\theta},be^{-i\theta},ce^{i\theta},ce^{-i\theta},de^{i\theta},de^{-i\theta};q)_\infty}
\frac{(ab;q)_\infty}{(ae^{i\theta},ae^{-i\theta};q)_\infty}{d}\theta\nonumber\\[6pt]
&=& \sum_{n=0}^{\infty}\int_{0}^{\pi}
\frac{(e^{2i\theta},e^{-2i\theta};q)_\infty}
{(be^{i\theta},be^{-i\theta},ce^{i\theta},ce^{-i\theta},de^{i\theta},de^{-i\theta};q)_\infty}
\frac{(f;q)_{n}g^n}{(q;q)_{n}}D_q^n\left\{\frac{(ab;q)_\infty}
{(ae^{i\theta},ae^{-i\theta};q)_\infty}\right\}{d}\theta.\nonumber
\end{eqnarray}
Let $U_n(\theta)$ denote the integrand in the last line of the above
equation. We make the assumption  $0 < q < 1$ so that, for $0 \leq
\theta \leq \pi$,
\begin{equation}
|(|x|;q)_{\infty}| \leq |(xe^{\pm i\theta};q)_{\infty}| \leq
(-|x|;q)_{\infty}
\end{equation}
and
\begin{equation}
|(e^{\pm 2i\theta};q)_\infty| \leq (-1;q)_{\infty}.
\end{equation}
Now we rewrite the series $\sum_{n=0}^{\infty}U_n(\theta)$ into
another form $\sum_{n=0}^{\infty}V_n(\theta)$ in order to prove its
uniform convergence. In the proof of Theorem \ref{gtef3}, one sees
that the
 absolute convergence of the ${}_3\phi_2$ series under the condition
 $|bs|,\,|b\,t|<1$
 implies the absolute convergence of the
sum
$$\sum_{n=0}^{\infty}\frac{(a;q)_nb^n}{(q;q)_n}
D_q^n\left\{\frac{(cv;q)_\infty}{(cs,ct;q)_\infty}\right\}.$$
Therefore, under the condition $|g|<1$, it follows from Theorem
\ref{gtef3} that
\begin{eqnarray}
\lefteqn{\sum_{n=0}^{\infty}\frac{(f;q)_{n}g^n}{(q;q)_{n}}
D_q^n\left\{\frac{(ab;q)_\infty}{(ae^{i\theta},
ae^{-i\theta};q)_\infty}\right\}}\nonumber\\[6pt]
&=&\frac{(fge^{i\theta},ab;q)_\infty}{(ge^{i\theta},
ae^{i\theta},ae^{-i\theta};q)_\infty}
{}_{3}\phi_2\left[\begin{array}{ccc}
f,&ae^{i\theta},&be^{i\theta}\\
&fge^{i\theta},&ab
\end{array};q,ge^{-i\theta}\right].
\end{eqnarray}
Hence
\begin{eqnarray}
\sum_{n=0}^{\infty}U_n(\theta)&=&\frac{(e^{2i\theta},e^{-2i\theta},fge^{i\theta},ab;q)_\infty}
{(ae^{i\theta},ae^{-i\theta},be^{i\theta},be^{-i\theta},ce^{i\theta},ce^{-i\theta},
de^{i\theta},de^{-i\theta},ge^{i\theta};q)_\infty}\nonumber\\[6pt]
&&\quad \times {}_{3}\phi_2\left[\begin{array}{ccc}
f,&ae^{i\theta},&be^{i\theta}\\
&fge^{i\theta},&ab
\end{array};q,ge^{-i\theta}\right]\nonumber\\[6pt]
&=&\frac{(e^{2i\theta},e^{-2i\theta},fge^{i\theta},ab;q)_\infty}
{(ae^{i\theta},ae^{-i\theta},be^{i\theta},be^{-i\theta},ce^{i\theta},ce^{-i\theta},
de^{i\theta},de^{-i\theta},ge^{i\theta};q)_\infty}\nonumber\\[6pt]
&&\quad \times\sum_{n=0}^{\infty}\frac{(f,ae^{i\theta},be^{i\theta};q)_n}
{(q,fge^{i\theta},ab;q)_n}
\left(ge^{-i\theta}\right)^n.
\end{eqnarray}
Now, let
\begin{eqnarray}
V_n(\theta)&=&\frac{(e^{2i\theta},e^{-2i\theta},fge^{i\theta},ab;q)_\infty}
{(ae^{i\theta},ae^{-i\theta},be^{i\theta},be^{-i\theta},ce^{i\theta},ce^{-i\theta},
de^{i\theta},de^{-i\theta},ge^{i\theta};q)_\infty}\nonumber\\[6pt]
&&\quad \times\frac{(f,ae^{i\theta},be^{i\theta};q)_n}{(q,fge^{i\theta},ab;q)_n}
\left(ge^{-i\theta}\right)^n.
\end{eqnarray}

 By the Weierstrass M-Test, it remains  to find a
convergent series $\sum_{n=0}^{\infty}M_n$, where $M_n$ is
independent of $\theta$, such that $|V_n(\theta)|\leq M_n$. For
$\max\{|a|,|b|,|c|,|d|,|g|\}<1$, we may choose
\begin{equation}
M_n=\left(\frac{(-1;q)_\infty}
{(|a|,|b|,|c|,|d|;q)_\infty}\right)^2\frac{(-|fg|,ab;q)_\infty}
{(|g|;q)_\infty}\frac{(-|f|,-|a|,-|b|;q)_{n}|g|^n}{|(q,|fg|,ab;q)_{n}|}.
\end{equation}
It is easy to see that $\sum_{n=0}^{\infty}M_n$ is convergent when
$|g|<1$. It follows that  the Cauchy operator commutes with the
integral in \eqref{44}, so \eqref{44} can be written as
\begin{eqnarray*}
\lefteqn{ \int_{0}^{\pi}\frac{(e^{2i\theta},e^{-2i\theta};q)_\infty}
{(be^{i\theta},be^{-i\theta},ce^{i\theta},ce^{-i\theta},de^{i\theta},de^{-i\theta};q)_\infty}
\sum_{n=0}^{\infty}\frac{(f;q)_{n}g^n}{(q;q)_{n}}
D_q^n\left\{\frac{(ab;q)_\infty}{(ae^{i\theta},ae^{-i\theta};q)_\infty}\right\}{d}\theta}
\nonumber \\[6pt]
&=& \int_{0}^{\pi}\frac{(e^{2i\theta},e^{-2i\theta};q)_\infty}
{(be^{i\theta},be^{-i\theta},ce^{i\theta},ce^{-i\theta},de^{i\theta},de^{-i\theta};q)_\infty}
T(f,g;D_q)\left\{\frac{(ab;q)_\infty}
{(ae^{i\theta},ae^{-i\theta};q)_\infty}\right\}{d}\theta.
\end{eqnarray*}
Finally, we may come to the general condition $|q|<1$ by the
argument of analytic continuation. Hence, under the condition
$\max\{|a|,|b|,|c|,|d|,|g|\}<1$, we have shown that it is valid to
exchange  the Cauchy operator and the integral when we apply the
Cauchy  operator to \eqref{integral}.

Now, applying $T(f,g;D_q)$ to \eqref{integral} with respect to the
parameter $a$ gives
\begin{eqnarray}
\lefteqn{\int_{0}^{\pi}\frac{(e^{2i\theta},e^{-2i\theta};q)_\infty}
{(be^{i\theta},be^{-i\theta},ce^{i\theta},ce^{-i\theta},
de^{i\theta},de^{-i\theta};q)_\infty}
\frac{(fge^{i\theta},ab;q)_\infty}{(ge^{i\theta},
ae^{i\theta},ae^{-i\theta};q)_\infty}}\nonumber\\[6pt]
&&\quad\quad\quad\times{}_{3}\phi_2\left[\begin{array}{ccc}
f,&ae^{i\theta},&be^{i\theta}\\[6pt]
&fge^{i\theta},&ab
\end{array};q,ge^{-i\theta}\right]{d}\theta\nonumber\\[6pt]
&=&\frac{2\pi}{(q,bc,bd,cd;q)_\infty}\frac{(cfg,abcd;q)_\infty}
{(cg,ac,ad;q)_\infty}
{}_{3}\phi_2\left[\begin{array}{ccc}
f,&ac,&bc\\
&cfg,&abcd
\end{array};q,dg\right],
\end{eqnarray}
where $\max\{|a|,|b|,|c|,|d|,|g|\}<1$. This implies the desired
formula. The proof is completed. \qed
\end{pf}

In fact, the above proof also implies the convergence of the
integral in Theorem \ref{ThmAskeyWilson}. Once it has been shown
that the sum commutes with the integral, one sees that the integral
obtained from exchanging the sum and the integral is convergent.

Setting $f=ab$ in \eqref{exAskeywilson}, by the $q$-Gauss sum
\cite[Appendix II.8]{Gasper-Rahman}:
\begin{equation}
{}_{2}\phi_1\left[\begin{array}{cc}
a,&b\\
&c
\end{array};q,\frac{c}{ab}\right]=\frac{(c/a,c/b;q)_\infty}
{(c,c/ab;q)_\infty},\ \ |c/ab|<1,\label{gauss}
\end{equation}
 we arrive at the Ismail-Stanton-Viennot integral
\eqref{ISV}.

Setting $f=abcd$ in \eqref{exAskeywilson},
by means of the $q$-Gauss sum \eqref{gauss}
we find the following formula which we have not seen in the
literature.

\begin{cor} We have
\begin{eqnarray}
\lefteqn{\int_{0}^{\pi}
\frac{(e^{2i\theta},e^{-2i\theta},abcdge^{i\theta};q)_\infty}
{(ae^{i\theta},ae^{-i\theta},be^{i\theta},
be^{-i\theta},ce^{i\theta},ce^{-i\theta},
de^{i\theta},de^{-i\theta},ge^{i\theta};q)_\infty}}\nonumber\\[6pt]
&&\quad\quad\times{}_{3}\phi_2\left[\begin{array}{ccc}
abcd,&ae^{i\theta},&be^{i\theta}\\
&abcdge^{i\theta},&ab
\end{array};q,ge^{-i\theta}\right]{d}\theta\nonumber\\[6pt]
&=&\frac{2\pi(abcd,acdg,bcdg;q)_\infty}{(q,ab,ac,ad,bc,bd,cd,cg,dg;q)_\infty},
\end{eqnarray}
where $\max\{|a|,|b|,|c|,|d|,|g|\}<1$.
\end{cor}

\section{A Further Extension of the Askey-Roy Integral}
\label{se Askey-Roy Integral}

Askey and Roy \cite{Askey-Roy} used Ramanujan's ${}_1\psi_1$
summation formula \cite[Appendix II.29]{Gasper-Rahman} to derive the
following integral formula:
\begin{eqnarray}\label{Askey-Roy}
\lefteqn{\frac{1}{2\pi}\int_{-\pi}^{\pi}\frac{(\rho
e^{i\theta}/d,qde^{-i\theta}/\rho,\rho
ce^{-i\theta},qe^{i\theta}/c\rho;q)_\infty}
{(ae^{i\theta},be^{i\theta},ce^{-i\theta},
de^{-i\theta};q)_\infty}{d}\theta \nonumber}\\[6pt]
&&=\frac{(abcd,\rho c/d,dq/\rho
c,\rho,q/\rho;q)_\infty}{(q,ac,ad,bc,bd;q)_\infty},
\end{eqnarray}
where $\max\{|a|,|b|,|c|,|d|\}<1$ and $cd\rho\neq0$, which is called
the Askey-Roy integral.

In \cite{Gasper}, Gasper discovered an integral formula
\begin{eqnarray}
\lefteqn{\frac{1}{2\pi}\int_{-\pi}^{\pi}\frac{(\rho
e^{i\theta}/d,qde^{-i\theta}/\rho,\rho
ce^{-i\theta},qe^{i\theta}/c\rho,abcdfe^{i\theta};q)_\infty}
{(ae^{i\theta},be^{i\theta},fe^{i\theta},ce^{-i\theta},de^{-i\theta};q)_\infty}
{d}\theta}\nonumber\\[6pt]
&=&\frac{(abcd,\rho c/d,dq/\rho
c,\rho,q/\rho,bcdf,acdf;q)_\infty}
{(q,ac,ad,bc,bd,cf,df;q)_\infty},\label{Gasper}
\end{eqnarray}
provided $\max\{|a|,|b|,|c|,|d|,|f|\}<1$ and $cd\rho\neq0$,  which
is an extension of the Askey-Roy integral. Note that Rahman and
Suslov \cite{RahmanSuslov} found a proof of Gasper's formula
\eqref{Gasper} based on the technique of iteration with respect to
the parameters of $\rho(s)$ in the integral
$$\int_{C}\rho(s)q^{-s}{d}s,$$
where $\rho(s)$ is the solution of a Pearson-type first-order
difference equation.

 In this section, we first derive an extension of the Askey-Roy
integral by applying the Cauchy  operator. We see that Gasper's
formula \eqref{Gasper} is a special case of this extension
\eqref{askey}. Moreover, a further extension of the Askey-Roy
integral can be obtained by taking the action of the Cauchy
 operator  on Gasper's formula.

\begin{thm}\label{t-Roy} We have
\begin{eqnarray}\label{Roy}
\lefteqn{\frac{1}{2\pi}\int_{-\pi}^{\pi}\frac{(\rho
e^{i\theta}/d,qde^{-i\theta}/\rho,\rho
ce^{-i\theta},qe^{i\theta}/c\rho,abcdfe^{i\theta},ghe^{i\theta};q)_\infty}
{(ae^{i\theta},be^{i\theta},fe^{i\theta},he^{i\theta},
ce^{-i\theta},de^{-i\theta};q)_\infty}}\nonumber\\[6pt]
&&\times{}_{3}\phi_2\left[\begin{array}{ccc}
g,&ae^{i\theta},&fe^{i\theta}\\
&ghe^{i\theta},&abcdfe^{i\theta}
\end{array};q,bcdh\right]  {d} \theta\nonumber\\[6pt]
&=&\frac{(abcd,\rho c/d,dq/\rho
c,\rho,q/\rho,bcdf,acdf,cgh;q)_\infty}
{(q,ac,ad,bc,bd,cf,ch,df;q)_\infty}\nonumber\\[6pt]
&&\quad \times{}_{3}\phi_2\left[\begin{array}{ccc}
g,&ac,&cf\\
&cgh,&acdf
\end{array};q,dh\right],
\end{eqnarray}
where $\max\{|a|,|b|,|c|,|d|,|f|,|h|\}<1$ and $cd\rho\neq0$.
\end{thm}

\begin{pf} As in the proof of the extension of the Askey-Wilson
integral, we can show that the Cauchy operator also commutes with
the Aksey-Roy integral. So we may apply the Cauchy operator
$T(f,g;D_q)$ to  both sides of the Askey-Roy integral
\eqref{Askey-Roy}  with respect to the parameter $a$. It follows
that
\begin{eqnarray}\label{askey}
\lefteqn{\frac{1}{2\pi}\int_{-\pi}^{\pi}\frac{(\rho
e^{i\theta}/d,qde^{-i\theta}/\rho,\rho
ce^{-i\theta},qe^{i\theta}/c\rho,fge^{i\theta};q)_\infty}
{(ae^{i\theta},be^{i\theta},ce^{-i\theta},de^{-i\theta},ge^{i\theta};q)_\infty}
{d}\theta}\nonumber\\[6pt]
&=&\frac{(abcd,cfg,\rho c/d,dq/\rho
c,\rho,q/\rho;q)_\infty}{(q,ac,ad,bc,bd,cg;q)_\infty}
{}_{3}\phi_2\left[\begin{array}{ccc}
f,&ac,&bc\\
&cfg,&abcd
\end{array};q,dg\right],
\end{eqnarray}
where $\max\{|a|,|b|,|c|,|d|,|g|\}<1$ and $cd\rho\neq0$.
\end{pf}

Putting $f= abcd$ and  $g= f$ in \eqref{askey}, by the $q$-Gauss sum
\eqref{gauss}, we get the formula  \eqref{Gasper} due to Gasper.

In order to apply the Cauchy  operator to Gasper's formula
\eqref{Gasper}, we rewrite it as
\begin{eqnarray}\label{Gasperch}
&&\frac{1}{2\pi}\int_{-\pi}^{\pi}\frac{(\rho
e^{i\theta}/d,qde^{-i\theta}/\rho,\rho
ce^{-i\theta},qe^{i\theta}/c\rho;q)_\infty}
{(be^{i\theta},fe^{i\theta},ce^{-i\theta},de^{-i\theta};q)_\infty}
\frac{(abcdfe^{i\theta};q)_\infty}{(ae^{i\theta},abcd;q)_\infty}
{d}\theta \nonumber\\[6pt]
&&\quad=\frac{(\rho c/d,dq/\rho
c,\rho,q/\rho,bcdf;q)_\infty}{(q,bc,bd,cf,df;q)_\infty}
\frac{(acdf;q)_\infty}{(ac,ad;q)_\infty}.
\end{eqnarray}
The proof is thus completed by employing the operator $T(g,h;D_q)$
with respect to the parameter $a$ to the above identity.\qed

Replacing $a$, $g$ by $g$, $cdfg$, respectively, and then taking
$h=a$ in \eqref{Roy}, we are led to the following identity due to
 Zhang and Wang \cite{Zhang-Wang}.

\begin{cor} We have
\begin{eqnarray}
&&\frac{1}{2\pi}\int_{-\pi}^{\pi}\frac{(\rho e^{i\theta}/d,
qde^{-i\theta}/\rho, \rho ce^{-i\theta}, qe^{i\theta}/c\rho,
abcdfge^{i\theta}, bcdfge^{i\theta};q)_{\infty}}{(ae^{i\theta},
be^{i\theta}, fe^{i\theta}, ge^{i\theta}, ce^{-i\theta},
de^{-i\theta};q)_{\infty}} \nonumber \\[6pt]
&&\quad\quad\quad\quad\times {}_{3}\phi_2\left[\begin{array}{ccc}
fe^{i\theta},&ge^{i\theta},&gcdf\\
&acdfge^{i\theta},&bcdfge^{i\theta}
\end{array};q,abcd\right]{d}\theta \nonumber \\[6pt]
&&\quad\quad\quad=\frac{(\rho c/d, dq/\rho c, \rho, q/\rho, acdf,
acdg, bcdf, bcdg, cdfg;q)_{\infty}}{(q, ac, ad, bc, bd, cf, df, cg,
dg;q)_{\infty}},
\end{eqnarray}
 where $\max\{|a|,|b|,|c|,|d|,|f|,|g|\}<1$ and $cd\rho\neq0$.
\end{cor}

\section{The Bivariate Rogers-Szeg\"o
Polynomials}\label{se Bivariate Rogers-Szego Polynomials}

In this section, we show that Mehler's formula and the Rogers
formula for the bivariate Rogers-Szeg\"o polynomials  can be easily
derived from the application of the Cauchy operator.
The bivariate Rogers-Szeg\"{o} polynomials are closely
related to the continuous big $q$-Hermite polynomials.
However, it seems that the following form of the
bivariate Rogers-Szeg\"o polynomials are
 introduced by Chen, Fu and Zhang \cite{ChenFuZhang}, as defined by
\begin{equation}
h_n(x,y|q)=\sum_{k=0}^{n}{n \brack k}P_k(x,y),
\end{equation}
where the Cauchy polynomials are given by
\[ P_k(x,y)=x^k(y/x;q)_k=(x-y)(x-qy)\cdots (x-q^{n-1}y),\]
which naturally arise in the $q$-umbral calculus. Setting $y=0$, the
polynomials $h_n(x,y|q)$ reduce to  the classical Rogers-Szeg\"{o}
polynomials $h_n(x|q)$ defined by
\begin{equation}
h_n(x|q)=\sum_{k=0}^{n}{n \brack k}x^k.
\end{equation}

It should be noted that Mehler's formula for the bivariate
Rogers-Szeg\"o polynomials is due to Askey, Rahman, and Suslov
\cite[Eq. (14.14)]{Askey-Rahman-Suslov}. They obtained the
nonsymmetric Poisson kernel formula for the continuous big
$q$-Hermite polynomials, often denoted by $H_n(x;a|q)$. The formula
of Askey, Rahman, and Suslov can be easily formulated in terms of
$h_n(x,y|q)$. Recently, Chen, Saad, and Sun presented an approach to
Mehler's formula and the Rogers formula for $h_n(x,y|q)$ by using
the homogeneous difference operator $D_{xy}$ introduced by Chen, Fu,
and Zhang. As will be seen, the Cauchy operator turns out to be more
efficient
 compared with the techniques used in \cite{Chen-Saad-Sun}.

We recall that the generating function of the bivariate
Rogers-Szeg\"{o} polynomials
\begin{equation}\label{gf}
\sum_{n=0}^{\infty}h_n(x,y|q)\frac{t^n}{(q;q)_n}=\frac{(yt;q)_\infty}
{(t,xt;q)_\infty},
\end{equation}
where $\max\{|x|,|xt|<1\}$, can be derived from the Euler identity
\eqref{Euler}
 using the Cauchy operator.

A direct calculation shows that
\begin{eqnarray}\label{DK2}
D_q^k\left\{a^n\right\} & = & \left\{\begin{array}{ll}
a^{n-k}(q^{n-k+1};q)_k,& 0\leq k\leq n,\\[6pt]
0, &k > n.\end{array}\right.
\end{eqnarray}

From the identity \eqref{DK2}, we can
 easily establish the following lemma.

\begin{lem} \label{m-l}
We have
\begin{equation}
T(a,b;D_q)\left\{c^n\right\}=\sum_{k=0}^{n}{n\brack
k}(a;q)_kb^kc^{n-k}.
\end{equation}
\end{lem}

Applying $T(a,b;D_q)$ to the Euler identity \eqref{Euler} with
respect to the parameter $z$, we get
\begin{equation}
\sum_{n=0}^{\infty}\frac{z^n}{(q;q)_n}\sum_{k=0}^{n}{n\brack
k}(a;q)_k\left(\frac{b}{z}\right)^k=\frac{(ab;q)_\infty}
{(b,z;q)_\infty},\label{rogers-szego}
\end{equation}
which leads to  \eqref{gf} by suitable substitutions.

The reason that we employ the Cauchy operator to deal with the
bivariate Rogers-Szeg\"o polynomials is based on the following fact
\begin{equation}
h_n(x,y|q)=\lim_{c \rightarrow 1}T(y/x,x;D_q)\left\{c^n\right\}.
\end{equation}

We are ready to describe how one can employ the
Cauchy operator to derive Mehler's formula and the Rogers
formula for $h_n(x,y|q)$.

\begin{thm}[Mehler's formula for $h_n(x,y|q)$] We have
\begin{equation}\label{gemehler}
\sum_{n=0}^{\infty}h_n(x,y|q)h_n(u,v|q)
\frac{t^n}{(q;q)_n}=\frac{(ty,tv;q)_\infty}{(t,tu,tx;q)_\infty}
{}_{3}\phi_2\left[\begin{array}{ccc}
t,&y/x,&v/u\\
&ty,&tv
\end{array};q,tux\right],
\end{equation}
where $\max\{|t|,|tu|,|tx|,|tux|\}<1$.
\end{thm}

\begin{pf}
By Lemma \ref{m-l}, the left hand side of (\ref{gemehler})
can be written as
\begin{eqnarray*}
\lefteqn{ \sum_{n=0}^{\infty}h_n(x,y|q)\lim_{c \rightarrow
1}T(v/u,u;D_q)\left\{c^n\right\}
\frac{t^n}{(q;q)_n}\nonumber}\\[6pt]
& = &
   \lim_{c \rightarrow 1}T(v/u,u;D_q)\left\{
\sum_{n=0}^{\infty}h_n(x,y|q)\frac{(ct)^n}{(q;q)_n}\right\}\nonumber.
\end{eqnarray*}
In view of the generating function (\ref{gf}), the above sum
equals
\begin{eqnarray}
\lefteqn{\lim_{c \rightarrow
1}T(v/u,u;D_q)\left\{\frac{(cty;q)_\infty}
{(ct,ctx;q)_\infty}\right\}\nonumber} \\[6pt]
&=&\lim_{c \rightarrow
1}\left(\frac{(tv,cty;q)_\infty}{(tu,ct,ctx;q)_\infty}
{}_{3}\phi_2\left[\begin{array}{ccc}
v/u,&ct,&y/x\\
&tv,&cty
\end{array};q,tux\right]\right)\nonumber\\[6pt]
&=&\frac{(ty,tv;q)_\infty}{(t,tu,tx;q)_\infty}
{}_{3}\phi_2\left[\begin{array}{ccc}
t,&y/x,&v/u\\
&ty,&tv
\end{array};q,tux\right],
\end{eqnarray}
where $\max\{|t|,|tu|,|xt|,|tux|\}<1$. This completes the proof.
\qed
\end{pf}

We see that \eqref{gemehler} is equivalent to \cite[Eq.
(2.1)]{Chen-Saad-Sun} in terms of  Sears' transformation formula
\eqref{hall}. Setting $y=0$ and $v=0$ in \eqref{gemehler} and
employing the Cauchy $q$-binomial theorem \eqref{Cauchy}, we obtain
Mehler's formula \cite{Chen1, IsmailStanton, Rogers1893,
stanton2000} for the Rogers-Szeg\"{o} polynomials.

\begin{cor}We have
\begin{equation}\label{mehler}
\sum_{n=0}^{\infty}h_n(x|q)h_n(u|q)
\frac{t^n}{(q;q)_n}=\frac{(t^2ux;q)_\infty}{(t,tu,tx,tux;q)_\infty},
\end{equation}
where $\max\{|t|,|tu|,|tx|,|tux|\}<1$.
\end{cor}

\begin{thm}[The Rogers
formula for $h_n(x,y|q)$]We have
\begin{equation}
\sum_{n=0}^{\infty}\sum_{m=0}^{\infty}h_{m+n}(x,y|q)
\frac{t^n}{(q;q)_n}\frac{s^m}{(q;q)_m}=\frac{(ty;q)_\infty}{(s,t,tx;q)_\infty}
{}_{2}\phi_1\left[\begin{array}{cc}
t,&y/x\\
&ty
\end{array};q,sx\right],\label{homo2}
\end{equation}
where $\max\{|s|,|t|,|sx|,|tx|\}<1$.
\end{thm}

\begin{pf} Using Lemma \ref{m-l},
the left hand side of \eqref{homo2} equals
\begin{eqnarray}
\lefteqn{\sum_{n=0}^{\infty}\sum_{m=0}^{\infty} \lim_{c \rightarrow
1}T(y/x,x;D_q)\left\{c^{m+n}\right\}\frac{t^n}{(q;q)_n}\frac{s^m}
{(q;q)_m}}\nonumber\\[6pt]
&=&\lim_{c \rightarrow
1}T(y/x,x;D_q)\left\{\sum_{n=0}^{\infty}\frac{(ct)^n}{(q;q)_n}
\sum_{m=0}^{\infty}\frac{(cs)^m}{(q;q)_m}\right\}
\end{eqnarray}
\begin{eqnarray}
&=&\lim_{c \rightarrow
1}T(y/x,x;D_q)\left\{\frac{1}{(cs,ct;q)_\infty}\right\}
\qquad\qquad\qquad\nonumber\\[6pt]
&=&\frac{(ty;q)_\infty}{(s,t,tx;q)_\infty}
{}_{2}\phi_1\left[\begin{array}{cc}
t,&y/x\\
&ty
\end{array};q,sx\right],
\end{eqnarray}
where $\max\{|s|,|t|,|sx|,|tx|\}<1$.\qed
\end{pf}

Note that \eqref{homo2} is equivalent to \cite[Eq.
(3.1)]{Chen-Saad-Sun} in terms of Heine's transformation formula
\cite[Appendix III.1]{Gasper-Rahman}. Setting $y=0$ in
\eqref{homo2}, by the Cauchy $q$-binomial theorem \eqref{Cauchy} we
get the Rogers formula \cite{Chen1, Rogers1893, Rogers1} for the
Rogers-Szeg\"{o} polynomials.

\begin{cor} We have
\begin{equation}\label{Rogers}
\sum_{n=0}^{\infty}\sum_{m=0}^{\infty}h_{m+n}(x|q)
\frac{t^n}{(q;q)_n}\frac{s^m}{(q;q)_m}=\frac{(stx;q)_\infty}
{(s,sx,t,tx;q)_\infty},
\end{equation}
where $\max\{|s|,|t|,|sx|,|tx|\}<1$.
\end{cor}

\section{An Extension of Sears' Formula}\label{se Sears' Formula}

In this section, we give an extension of  the Sears two-term
summation formula \cite[Eq. (2.10.18)]{Gasper-Rahman}:
\begin{eqnarray}\label{sears}
\lefteqn{\int_{c}^{d}\frac{(qt/c,qt/d,abcdet;q)_\infty}
{(at,bt,et;q)_\infty}{d}_qt\nonumber}\\
&=&\frac{d(1-q)(q,dq/c,c/d,abcd,bcde,acde;q)_\infty}
{(ac,ad,bc,bd,ce,de;q)_\infty},
\end{eqnarray}
where $\max\{|ce|,|de|\}<1$.

From the Cauchy operator, we deduce the following extension of
(\ref{sears}).

\begin{thm} We have
\begin{eqnarray}\label{tSears}
\lefteqn{\int_{c}^{d}\frac{(qt/c,qt/d,abcdet,fgt;q)_\infty}
{(at,bt,et,gt;q)_\infty} {}_{3}\phi_2\left[\begin{array}{ccc}
f,&at,&et\\
&fgt,&abcdet
\end{array};q,bcdg\right] {d}_qt}\nonumber\\[6pt]
&=&\frac{d(1-q)(q,dq/c,c/d,abcd,bcde,acde,cfg;q)_\infty}
{(ac,ad,bc,bd,ce,cg,de;q)_\infty}\nonumber\\[6pt]
&&\times{}_{3}\phi_2\left[\begin{array}{ccc}
f,&ac,&ce\\
&cfg,&acde
\end{array};q,dg\right],
\end{eqnarray}
where $\max\{|bcdg|,|ce|,|cg|,|de|,|dg|\}<1$.
\end{thm}

\begin{pf}
We may rewrite \eqref{sears}  as
\begin{eqnarray}
\lefteqn{\int_{c}^{d}\frac{(qt/c,qt/d;q)_\infty}
{(bt,et;q)_\infty}\frac{(abcdet;q)_\infty}{(at,abcd;q)_\infty}
\text{d}_qt}\nonumber\\[6pt]
&=&\frac{d(1-q)(q,dq/c,c/d,bcde;q)_\infty}
{(bc,bd,ce,de;q)_\infty}\frac{(acde;q)_\infty}{(ac,ad;q)_\infty}.
\end{eqnarray}
Applying the operator $T(f,g;D_q)$  with respect to the parameter
$a$, we obtain \eqref{tSears}.\qed
\end{pf}

 As far as the convergence is concerned, the above integral is of
the following form
\begin{equation}\label{abn}
\sum_{n=0}^{\infty}A(n)\sum_{k=0}^{\infty}B(n,k).
\end{equation}
To ensure that the series (\ref{abn}) converges absolutely, we
assume that the following two conditions are satisfied:
\begin{enumerate}
\item
$\sum_{k=0}^{\infty}B(n,k)$ converges to $C(n)$, and  $C(n)$
   has a nonzero limit as
$n\rightarrow \infty$.
\item $\lim\limits_{n \rightarrow
\infty}|\frac{A(n)}{A(n-1)}|<1$.
\end{enumerate}

It is easy to see that under the above assumptions, (\ref{abn})
converges absolutely, since
\[
\lim_{n \rightarrow
\infty}\left|\frac{A(n)C(n)}{A(n-1)C(n-1)}\right| = \lim_{n
\rightarrow \infty}\left|\frac{A(n)}{A(n-1)}\right|<1 .
\]

It is easy to verify the double summations in \eqref{tSears} satisfy
the two assumptions of \eqref{abn}, so the convergence is
guaranteed.

\section{Extensions of $q$-Barnes' Lemmas}\label{se q-Barnes' Lemmas}

 In this section, we obtain extensions of the $q$-analogues of
Barnes' lemmas. Barnes' first lemma \cite{Barne1} is an integral
analogue of Gauss'
 ${}_2F_1$ summation formula.
 Askey and Roy \cite{Askey-Roy} pointed out that Barnes' first
lemma is also an extension of the beta integral. Meanwhile,
Barnes'
second lemma \cite{Barne2} is an integral analogue of
Saalsch\"{u}tz's formula.

The following
 $q$-analogue of Barnes' first lemma is due to Watson, see \cite[Eq.
(4.4.3)]{Gasper-Rahman}:
\begin{eqnarray}\label{barnes1}
&&\frac{1}{2\pi i}\int_{-i\infty}^{i\infty}
\frac{(q^{1-c+s},q^{1-d+s};q)_\infty}{(q^{a+s},q^{b+s};q)_\infty}
\frac{\pi q^s {d}s}{\sin \pi(c-s) \sin\pi(d-s)} \nonumber \\[6pt]
&&\quad\quad= \frac{q^c}{\sin \pi(c-d)}
\frac{(q,q^{1+c-d},q^{d-c},q^{a+b+c+d};q)_{\infty}}
{(q^{a+c},q^{a+d},q^{b+c},q^{b+d};q)_{\infty}}.
\end{eqnarray}
The $q$-analogue of Barnes' second lemma is due to Agarwal, see
\cite{Agarwal} and  \cite[Eq. (4.4.6)]{Gasper-Rahman}:
\begin{eqnarray}\label{barnes2}
&&\frac{1}{2\pi i}\int_{-i\infty}^{i\infty}
\frac{(q^{1+s},q^{d+s},q^{1+a+b+c+s-d};q)_\infty}
{(q^{a+s},q^{b+s},q^{c+s};q)_\infty}\
\frac{\pi q^s {d}s}{\sin \pi s \sin\pi(d+s)} \nonumber \\[6pt]
&&\quad\quad=\csc \pi d\
\frac{(q,q^d,q^{1-d},q^{1+b+c-d},q^{1+a+c-d},q^{1+a+b-d};q)_{\infty}}
{(q^a,q^b,q^c,q^{1+a-d},
q^{1+b-d},q^{1+c-d};q)_{\infty}},
\end{eqnarray}
where $\text{Re}\{s\log q-\log(\sin\pi s\sin\pi(d+s))\}<0$ for large
$|s|$. Throughout this section, the contour of integration always ranges from
$-i\infty$ to $i\infty$ so that the increasing sequences of poles of
integrand lie to the right and the decreasing sequences of poles
lie to the left of the contour, see \cite[p. 119]{Gasper-Rahman}. In order to ensure that the Cauchy
operator commutes with the integral,  we assume that
$q=e^{-\omega},\ \omega>0$.

We obtain the following extension of Watson's $q$-analogue of
Barnes' first lemma.

\begin{thm} We have
\begin{eqnarray}\label{barnes1-1}
\lefteqn{\frac{1}{2\pi i}\int_{-i\infty}^{i\infty}
\frac{(q^{1-c+s},q^{1-d+s},q^{e+f+s};q)_\infty}
{(q^{a+s},q^{b+s},q^{f+s};q)_\infty}
\frac{\pi q^s {d}s}{\sin \pi(c-s) \sin\pi(d-s)}} \nonumber \\[6pt]
&=&\frac{q^c}{\sin \pi(c-d)}
\frac{(q,q^{1+c-d},q^{d-c},q^{a+b+c+d},q^{c+e+f};q)_{\infty}}
{(q^{a+c},q^{a+d},q^{b+c},q^{b+d},q^{c+f};q)_{\infty}}\nonumber\\[6pt]
&&\times{}_{3}\phi_2\left[\begin{array}{ccc}
q^e,&q^{a+c},&q^{b+c}\\
&q^{c+e+f},&q^{a+b+c+d}
\end{array};q,q^{d+f}\right],
\end{eqnarray}
where $\max\{|q^f|,|q^{c+f}|,|q^{d+f}|\}<1$.
\end{thm}
\begin{pf}
Applying the operator $T(q^e,q^f;D_q)$ to \eqref{barnes1} with
respect to the parameter $q^a$, we arrive at \eqref{barnes1-1}. \qed
\end{pf}

Let us consider the special case when $e=a+b+c+d$.
 The ${}_{3}\phi_2$ sum on the right hand side of
\eqref{barnes1-1} turns out to be a ${}_{2}\phi_1$ sum and can be
summed by the $q$-Gauss formula \eqref{gauss}. Hence we get the
following formula derived by Liu \cite{Liu}, which is also an
extension of $q$-Barnes' first Lemma.

\begin{cor}\label{barnes-cor1} We have
\begin{eqnarray}\label{barnes1-3}
\lefteqn{\frac{1}{2\pi i}\int_{-i\infty}^{i\infty}
\frac{(q^{1-c+s},q^{1-d+s},q^{a+b+c+d+f+s};q)_\infty}
{(q^{a+s},q^{b+s},q^{f+s};q)_\infty}
\frac{\pi q^s {d}s}{\sin \pi(c-s) \sin\pi(d-s)}} \nonumber \\[6pt]
&=&\frac{q^c}{\sin \pi(c-d)}
\frac{(q,q^{1+c-d},q^{d-c},q^{a+b+c+d},q^{a+c+d+f},q^{b+c+d+f};q)_{\infty}}
{(q^{a+c},q^{a+d},q^{b+c},q^{b+d},q^{c+f},q^{d+f};q)_{\infty}},
\end{eqnarray}
where $\max\{|q^f|,|q^{c+f}|,|q^{d+f}|\}<1$.
\end{cor}

Clearly, \eqref{barnes1-3} becomes $q$-Barnes' first Lemma
\eqref{barnes1} for $f\rightarrow \infty$.
Based on  Corollary \ref{barnes-cor1}, employing the Cauchy
 operator again, we derive the following further
extension of  $q$-Barnes' first Lemma.

\begin{thm} We have
\begin{eqnarray}\label{barnes1-4}
\lefteqn{\frac{1}{2\pi i}\int_{-i\infty}^{i\infty}
\frac{(q^{1-c+s},q^{1-d+s},q^{a+b+c+d+f+s},q^{e+g+s};q)_\infty}
{(q^{a+s},q^{b+s},q^{f+s},q^{g+s};q)_\infty}
\frac{\pi q^s }{\sin \pi(c-s) \sin\pi(d-s)}} \nonumber \\[6pt]
&&\quad\times{}_{3}\phi_2\left[\begin{array}{ccc}
q^e,&q^{a+s},&q^{b+s}\\
&q^{e+g+s},&q^{a+b+c+d+f+s}
\end{array};q,q^{c+d+f+g}\right]{d}s\nonumber\\[6pt]
&=&\frac{q^c}{\sin \pi(c-d)}
\frac{(q,q^{1+c-d},q^{d-c},q^{a+b+c+d},q^{a+c+d+f},
q^{b+c+d+f},q^{c+e+g};q)_{\infty}}
{(q^{a+c},q^{a+d},q^{b+c},q^{b+d},q^{c+f},
q^{c+g},q^{d+f};q)_{\infty}}\nonumber\\[6pt]
&&\quad\quad\times{}_{3}\phi_2\left[\begin{array}{ccc}
q^e,&q^{a+c},&q^{b+c}\\
&q^{c+e+g},&q^{a+b+c+d}
\end{array};q,q^{d+g}\right],
\end{eqnarray}
where
$\max\{|q^f|,|q^g|,|q^{c+f}|,|q^{c+g}|,|q^{d+f}|,|q^{d+g}|,|q^{c+d+f+g}|\}<1$.
\end{thm}

 We conclude this paper with the following extension of
 Agarwal's $q$-analogue of Barnes' second lemma. The proof
 is omitted.

\begin{thm} We have
\begin{eqnarray}\label{barnes2-1}
\lefteqn{\frac{1}{2\pi i}\int_{-i\infty}^{i\infty}
\frac{(q^{1+s},q^{d+s},q^{1+a+b+c+s-d},q^{e+f+s};q)_\infty}
{(q^{a+s},q^{b+s},q^{c+s},q^{f+s};q)_\infty}\
\frac{\pi q^s }{\sin \pi s \sin\pi(d+s)} }\nonumber \\[6pt]
&&\times{}_{3}\phi_2\left[\begin{array}{ccc}
q^e,&q^{a+s},&q^{b+s}\\
&q^{e+f+s},&q^{1+a+b+c+s-d}
\end{array};q,q^{1+c+f-d}\right]{d}s\nonumber\\[6pt]
&=&\csc \pi d\
\frac{(q,q^d,q^{1-d},q^{1+b+c-d},q^{1+a+c-d},q^{1+a+b-d},q^{e+f};q)_{\infty}}
{(q^a,q^b,q^c,q^{f},q^{1+a-d},
q^{1+b-d},q^{1+c-d};q)_{\infty}}\nonumber\\[6pt]
&&\times{}_{3}\phi_2\left[\begin{array}{ccc}
q^a,&q^{b},&q^{e}\\
&q^{e+f},&q^{1+a+b-d}
\end{array};q,q^{1+f-d}\right],
\end{eqnarray}
where $\max\{|q^f|,|q^{1+f-d}|,|q^{1+c+f-d}|\}<1$ and
$\text{Re}\{s\log q-\log(\sin\pi s\sin\pi(d+s))\}<0$ for large
$|s|$.
\end{thm}

\vspace{.2cm} \noindent{\bf Acknowledgments.}  We would like to
thank the referee and Lisa H. Sun for helpful comments. This work
was supported by the 973 Project, the PCSIRT Project of the Ministry
of Education, the Ministry of Science and Technology and the
National Science Foundation of China.

\end{document}